\documentclass[
11pt]{amsart}
\usepackage{amsthm, amsfonts, amssymb}

\theoremstyle{definition}
\newtheorem{ntn}{Notation}[section]
\newtheorem{dfn}[ntn]{Definition}
\theoremstyle{plain}
\newtheorem{lem}[ntn]{Lemma}
\newtheorem{prp}[ntn]{Proposition}
\newtheorem{thm}[ntn]{Theorem}
\newtheorem{cor}[ntn]{Corollary}
\theoremstyle{remark}
\newtheorem{rem}[ntn]{Remark}

\def\floor[#1]{\lfloor #1 \rfloor }

\newcommand{\z}{\mathbb{Z}}

\newcommand{\lan}{\langle}
\newcommand{\ran}{\rangle}
\newcommand{\fs}{F\langle S\rangle}
\newcommand{\htf}{{\rm ht}_{F^{op}}(v)}
\newcommand{\op}{{op}}
\newcommand{\fop}{F^{op}}
\newcommand{\sr}{{\rm sr}(R)}
\newcommand{\hht}{{\rm ht}}
\newcommand{\rn}{R^{2n}}
\newcommand{\GL}{\mathit{GL}}
\newcommand{\Spp}{\mathit{Sp}}
\newcommand{\ESp}{\mathit{ESp}}
\newcommand{\iur}{\mathcal{IU}(R^{2n})}
\newcommand{\mur}{\mathcal{MU}(R^{2n})}
\newcommand{\hur}{\mathcal{HU}(R^{2n})}
\newcommand{\iu}{\mathcal{IU}}
\newcommand{\hu}{\mathcal{HU}}
\newcommand{\n}{\mathcal{N}}

\renewcommand{\a}{\mathcal{A}}

\newcommand{\lii}{\mathit{Link}^{+}_X(x)}

\renewcommand{\o}{\mathcal{O}}
\renewcommand{\H}{\tilde{H}}
\newcommand{\ee}{\mathcal{E}}
\newcommand{\ff}{\mathcal{F}}
\renewcommand{\gg}{\mathcal{G}}

\renewcommand{\lll}{\mathcal{L}}
\renewcommand{\u}{\mathcal{U}}

\newcommand{\s}{\Sigma}
\newcommand{\si}{\sigma}

\newcommand{\dl}{\Delta}

\newcommand{\arr}{\rightarrow}
\newcommand{\harr}{\hookrightarrow}

\newcommand{\se}{\subseteq}
\newcommand{\bs}{\backslash}
\newcommand{\bcu}{\bigcup}

\newcommand{\mt}{\mapsto}
\newcommand{\two}{\twoheadrightarrow}

\DeclareMathAlphabet{\mathds}{U}{dsrom}{m}{n}
\DeclareMathAlphabet{\mathsc}{U}{rsfs}{m}{n}

\begin{document}

\title{Homology stability for symplectic groups}
\author{B. Mirzaii}
\author{W. van der Kallen}

\begin{abstract}
In this paper the homology stability for symplectic
groups over a ring with finite stable rank is established.
First we develop a `nerve theorem'
 on the homotopy type of a poset in terms
of a cover by subposets, where the cover is itself indexed by a poset.
We use the nerve theorem
 to show that a poset of sequences of isotropic vectors is
highly connected, as conjectured by Charney in the eighties.
\end{abstract}

\maketitle
\maketitle

\section{Introduction}

Interest in homological stability problems in algebraic $K$-theory
started with Quillen, who used it in \cite{quil2} to study the
higher $K$-groups of a
ring of integers. As a result of stability he proved that
these groups are finitely
generated (see also \cite{gra}). After that there has been considerable
interest in homological stability for general linear groups. The most
general results in this direction are due to the second author \cite{kall}
and Suslin \cite{sus}.

Parallel to this, similar questions for other classical groups
such as orthogonal and symplectic groups have been studied.
For work in this direction, see \cite{vog}, \cite{bet0}, \cite{bet},
\cite{ch}, \cite{pan1}, \cite{pan2}. The most general result
is due to Charney \cite{ch}. She proved the homology
stability for orthogonal and symplectic groups over a Dedekind domain.
Panin in \cite{pan2} proved a similar result but with a different
method and with better range of stability. There he gave a complete
formulation of the homology stability conjecture for these groups. Charney
in \cite{ch} also showed that  knowing the higher connectivity of
the poset of isotropic unimodular sequences $\iur$, with respect to the
given bilinear form on $\rn$,
one can get the homology stability of the corresponding family of
classical groups.

Our goal in this paper is to prove that homology stabilizes of
the symplectic
groups over rings with finite stable rank. To do so we
prove that the poset of isotropic unimodular sequences is indeed
highly connected. Recall that
Panin in \cite{pan1} had already sketched how one can do
this for a finite dimensional affine algebra over an infinite field.
However, while the assumption about the infinite field provides
a significant simplification,
it excludes
cases of primary interest, namely rings that are finitely generated
over the integers.

Our approach is as follows.
We first extend a
theorem of Quillen \cite[Thm 9.1]{quil7} which was his main tool
to prove that certain posets are highly connected.
We use it to develop a quantitative analogue for posets of the nerve
theorem, which expresses the homotopy type of a space in terms of the
the nerve of a suitable cover.
In our situation both the elements of the cover and the nerve are replaced
with posets.
We work with posets of ordered sequences `satisfying the chain condition',
as this is a good replacement for simplicial complexes in the presence
of group actions. (This roughly corresponds with
 the passage to a barycentric subdivision of a simplicial complex.)
The new nerve theorem allows us to
exploit the higher connectivity of the poset
of unimodular sequences due to the second author. The
higher connectivity of
the poset of isotropic unimodular sequences
follows inductively. We conclude
with the homology stability theorem.

\section{Preliminaries}

 Recall that a topological space $X$ is $(-1)$-{\it connected}
if it is non-empty, $0$-{\it connected} if it is non-empty and
path connected, $1$-{\it connected} if it is non-empty and simply
connected.  In general for $n \ge 1$, $X$ is called
$n$-{\it connected} if $X$ is nonempty, $X$ is $0$-connected and
$\pi_i(X, x)=0$ for every base point $x \in X$ and $1 \le i \le n$.
For $n \ge -1$ a space $X$ is called $n$-{\it acyclic} if it is
nonempty and $\H_i(X, \z)=0$ for $0 \le i \le n$. For $n<-1$ the
conditions of $n$-connectedness and $n$-acyclicness are vacuous.

\begin{thm}[Hurewicz]\label{hur}
For $n \ge 0 $, a topological space $X$ is
$n$-connected if and only if the reduced homology groups
$\H_i(X, \z)$ are trivial for $0 \le i \le n$ and $X$ is
$1$-connected if $n \ge 1$.
\end{thm}
\begin{proof}
See \cite{whi}, Chap. IV, Corollaries 7.7 and 7.8.
\end{proof}

Let $X$ be a partially ordered set or briefly a {\it poset}.
Consider the simplicial complex associated to $X$,
that is the simplicial complex where vertices or $0$-simplices
are the elements of $X$ and the $k$-simplices are the $(k+1)$-tuples
$(x_0, \dots, x_k)$ of elements of  $X$ with $x_0< \dots < x_k$.
We denote it again by $X$. We denote the geometric realization of
$X$ by $|X|$ and we consider it with the weak topology. It is
well known that $|X|$ is a CW-complex \cite{mil}. By a {\it morphism}
or {\it map} of posets $f: X \arr Y$ we mean an order-preserving
map i.~e.\ if $x \le x'$ then $f(x) \le f(x')$. Such a map induces
a continuous map $|f|:|X| \arr |Y|$.

\begin{rem}
If $K$ is a simplicial complex and $X$
the ordered set of simplices of $K$, then the space
$|X|$ is the barycentric subdivision of $K$. Thus every simplicial
complex, with weak topology, is homeomorphic to the geometric
realization of some, and in fact many, posets. Furthermore
since it is well known that any CW-complex is homotopy equivalent
to a simplicial complex, it follows that any interesting homotopy
type is realized as the geometric realization of a poset.
\end{rem}

\begin{prp}\label{qui}
Let $X$ and $Y$ be posets.
\par {\rm (i)} $(${\rm Segal \cite{seg}}$)$ If $f, g : X \arr Y$ are maps
of posets such that $f(x) \le g(x)$ for all $x \in X$, then $|f|$ and
$|g|$ are homotopic.
\par {\rm (ii)} If the poset X has a minimal or maximal element then
$|X|$ is contractible.
\par {\rm (iii)} If $X^{op}$ denotes the opposite poset of X, i.~e.\
with opposite ordering, then $|X^{op}| \simeq |X|$.
\end{prp}
\begin{proof}
\par (i) Consider the poset $I=\{0,1: 0<1\}$ and define the poset map
$h: I \times X \arr Y$ as $h(0, x)=f(x)$, $h(1, x)=g(x)$. Since $|I|\simeq
[0,1]$, we have $|h| :[0,1] \times |X| \arr |Y|$ with
$|h|(0, x)=|f|(x)$ and  $|h|(1, x)=|g|(x)$. This shows that $|f|$
and $|g|$ are homotopic.
\par (ii) Suppose $X$ has a maximal element $z$. Consider the map
$f: X \arr X$ with $f(x)=z$ for every $x \in X$. Clearly for
every $x \in X$, ${\rm id}_X(x) \le f(x)$. This shows that
${\rm id}_X$ and the constant map $f$ are homotopic. So $X$ is
contractible. If $X$ has a minimal element the proof is similar.
\par (iii). This is natural and easy.
\end{proof}

The construction $X\mapsto |X|$ allows us to assign topological
concepts to posets. For example we define the homology groups
of a poset $X$ to be those of $|X|$, we call $X$ $n$-{\it connected}
or {\it contractible} if $|X|$ is $n$-connected or contractible
etc. Note that $X$ is connected if and only if $X$ is connected as
a poset. By the {\it dimension} of a poset $X$,
we mean the dimension of the
space $|X|$, or equivalently the supremum of the
integers $n$ such that there is a chain $x_0 < \cdots <x_n$ in
$X$. By convention the empty set has dimension $-1$.

Let $X$ be  a poset and $x \in X$. Define
$Link^+_X(x) := \{ u \in X : u > x\}$ and
$Link^-_X(x) := \{ u \in X : u < x\}$.
Given a map $f: X \arr Y$ of posets and an element $y \in Y$,
define subposets $f/y$ and $y\bs f$ of $X$ as follows
\begin{gather*}
f/y:=\{ x \in X : f(x) \le y\} \qquad
y\bs f := \{ x \in X : f(x) \ge y\}.
\end{gather*}
In fact $f/y=f^{-1}(Y_{\le y})$ and
$y\bs f=f^{-1}(Y_{\ge y})$ where $Y_{\le y}=\{ z \in Y: z \le y\}$
and $Y_{\ge y}=\{ z \in Y: z \ge y\}$. Note that by \ref{qui} (ii),
$Y_{\le y}$ and $Y_{\ge y}$ are contractible. If
${\rm id}_Y: Y \arr Y$ is the identity map, then ${\rm id}_Y/y=Y_{\le y}$ and
${\rm id}_Y\bs y=Y_{\ge y}$.

Let $\ff: X \arr  \underline{{\rm Ab}}$ be a functor from a poset
$X$, regarded as a category in usual way,
to the category of abelian groups. We define the homology groups
$H_i(X, \ff)$ of $X$ with coefficient $\ff$ to be the homology of the
complex $C_\bullet(X, \ff)$ given by
\[
C_n(X, \ff) =\underset{x_0<\cdots<x_n}{\bigoplus}\ff(x_0)
\]
where the direct sum is taken over all $n$-simplices in $X$, with
differential $\partial_n=\s_{i=0}^n(-1)^id^n_i$ where
$d^n_i: C_n(X, \ff) \arr C_{n-1}(X, \ff)$ and $d^n_i$ takes
the $(x_0<\cdots<x_n)$-component of $C_n(X, \ff)$ to the
$(x_0< \cdots < \widehat{x_i}< \cdots <x_n)$-component of
$C_{n-1}(X, \ff)$ via $d^n_i = id_{\ff(x_0)}$ if $i>0$ and
$d^n_0 :\ff(x_0)\arr \ff(x_1)$. In particular, for the empty set we have
$H_i(\varnothing, \ff)=0$ for $i \ge 0$.

Let $\ff$ be the constant functor $\z$. Then the homology groups
with this coefficient coincide with the integral homology of $|X|$,
that is $H_k(X, \z)= H_k(|X|, \z)$ for all $k \in \z$, \cite[App. II]{g-z}.
Let $\H_i(X, \z)$ denote the reduced integral
homology of the poset $X$, that is
$\tilde{H}_i(X, \z)={\rm ker}\{H_i(X, \z) \arr H_i(pt, \z)\}$
if $X\neq \varnothing$ and
$\tilde{H}_i(\varnothing, \z)= \begin{cases}\z & \text{if $i=-1$}\\
0 & \text{if $i\neq -1$} \end{cases}$.
So $\H_i(X, \z)=H_i(X, \z)$ for $i \ge 1$ and for $i=0$ we have the
exact sequence
\[
0 \arr \H_0(X, \z) \arr H_0(X, \z) \arr \z
\arr \tilde{H}_{-1}(X, \z) \arr 0
\]
where $\z$ is identified with the group $H_0(pt, \z)$. Notice
that $H_0(X, \z)$ is identified with the free abelian group
generated by the connected components of $X$.

A {\it local system} of abelian groups on a space (resp. poset) $X$ is a
functor $\ff$ from the groupoid of $X$ (resp.\ viewed  $X$ as a category),
to the category of abelian groups which is morphism-inverting, i.~e.\
such that the map $\ff(x) \arr \ff(x')$ associated to a path from $x$
to $x'$ (resp. $x \le x'$) is an isomorphism.
Clearly, a local system $\ff$ on a path connected space (resp. 0-connected
poset) is determined, up to canonical isomorphism, by the following 
data: if $x \in X$ is a base point, it suffices
to be given the group $\ff(x)$ and an action of $\pi_1(X, x)$ on $\ff(x)$.

The homology groups $H_k(X, \ff)$ of a space with a local system $\ff$
are a generalization
of the ordinary homology groups \cite[Chap. VI]{whi}.  In fact if
$X$ is a 0-connected space and if $\ff$ is a constant local system
on X, then $H_k(X, \ff) \simeq H_k(X, \ff(x_0))$ for every $x_0 \in X$
\cite[Chap. VI, 2.1]{whi}.

Let $X$ be a poset and $\ff$ a local system on $|X|$. Then the restriction
of $\ff$ to $X$ is a local system on $X$. Considering $\ff$
as a functor from $X$ to the category of abelian groups, we can
define $H_k(X, \ff)$ as in the above. Conversely if $\ff$ is a local
system on the poset $X$, then there is a unique 
local system, up to isomorphism, on $|X|$ such that
the restriction to $X$ is $\ff$ \cite[ Chap. VI, Thm 1.12]{whi},
\cite[I, Prop. 1]{quil3}. We denote both local systems by $\ff$.

\begin{thm}\label{q-w}
Let $X$ be a poset and $\ff$ a local system on $X$.
Then the homology groups $H_k(|X|, \ff)$ are
isomorphic with the homology groups $H_k(X, \ff)$.
\end{thm}
\begin{proof}
See \cite[Chap. VI, Thm. 4.8]{whi} or \cite[I, p.~91]{quil3}.
\end{proof}

\begin{thm}\label{wh1}
Let X be a path connected space with a base point $x$ and let
$\ff$ be a local system on X. Then the inclusion $\{x\} \harr X$
induces an isomorphism
$\ff(x)/G \overset{\simeq}{\longrightarrow} H_0(X, \ff)$
where G is the subgroup of $\ff(x)$ generated by all the elements of
the form $a- \beta a$ with $ a \in \ff(x)$, $\beta \in \pi_1(X, x)$.
\end{thm}
\begin{proof}
See \cite{whi}, Chap. VI, Thm. 2.8$^\ast$ and Thm. 3.2.
\end{proof}

We need the following interesting and well known lemma about the covering
spaces of the space $|X|$, where $X$ is a poset (or
more generally a simplicial set).
For a definition of a covering space, useful for our purpose, and some
more information, see \cite[Chap. 2]{spa}.

\begin{lem}\label{g-z-q}
The category of the covering spaces of the space $|X|$ of a poset $X$ is
equivalent to the category $\lll_S(X)$, the category of functors
$\ff: X \arr \underline{{\rm Set}}$, where $\underline{{\rm Set}}$ is
the category of sets, such that $\ff(x) \arr \ff(x')$ is a bijection
for every relation $x \le x'$.
\end{lem}
\begin{proof}
See \cite[I, p.~90]{quil3}. For the same proof with more
details see \cite[lem. 6.1]{sri}.
\end{proof}

\section{Homology and homotopy of posets}\label{hompos}

\begin{thm}\label{g-z}
Let $f: X \arr Y$ be a map of posets. Then there is a first
quadrant spectral sequence
\begin{gather*}
E_{p, q}^2 = H_p(Y, y \mt H_q(f/y, \z)) \Rightarrow H_{p+q}(X, \z).
\end{gather*}
The spectral sequence is functorial, in the sense that if
there is a commutative diagram of posets
\begin{gather*}
\begin{array}{ccc}
X'  & \overset{f'}{\longrightarrow} & Y'    \\
\Big\downarrow\vcenter{%
\rlap{$\scriptstyle{g_X}$}} &      &
\Big\downarrow\vcenter{%
\rlap{$\scriptstyle{g_Y}$}}\\
X          & \overset{f}{\longrightarrow}  &  Y
\end{array}
\end{gather*}
then there is a natural map from the spectral sequence arising from $f'$
to the spectral sequence arising from $f$. Moreover the  map
${g_X}_\ast: H_\ast(X', \z) \arr H_\ast(X, \z)$ is compatible with
this natural map.
\end{thm}
\begin{proof}
Let $C_{\ast, \ast}(f)$ be the double complex
such that $C_{p, q}(f)$ is the free abelian group generated by the set
$\{(x_0< \dots <x_q, f(x_q)<y_0< \dots < y_p): x_i \in X, y_i \in Y\}$.
The first spectral sequences of this double complex has as $E^1$-term
$E_{p, q}^1({\rm I})= \bigoplus_{y_0< \dots <y_p}H_q(f/y_0, \z)$.
By the general theory of double complexes (see for example
\cite[Chap.~5]{weib}),
we know that $E_{p, q}^2({\rm I})$ is the homology of the
chain complex $C_\ast(Y, \gg_q)=E_{\ast, q}^1({\rm I})$ where
$\gg_q: Y \arr \underline{{\rm Ab}}$, $\gg_q(y)=H_q(f/y, \z)$ and
hence $E_{p, q}^2({\rm I})= H_p(Y, y \mt H_q(f/y, \z))$.
The second spectral sequence has as $E^1$-term
$E_{p, q}^1({\rm II})= \bigoplus_{f(x_q)<y_0< \dots < y_p}
H_p({\rm id}_Y \bs f(x_q), \z)$.
But by \ref{qui} (ii), ${\rm id}_Y \bs f(x_q)=f(x_q)\bs Y$ is contractible, so
$E_{\ast, 0}^1({\rm II})=C_\ast(X^\op, \z)$ and $E_{\ast, q}^1({\rm II})=0$
for $q > 0$. Hence
$H_\ast(C_{\ast, \ast}(f))\simeq H_\ast(X^\op, \z) \simeq  H_\ast(X, \z)$.
This completes the proof of existence and convergence of the spectral sequence.
The functorial behavior of the spectral sequence follows from the functorial
behavior of the spectral sequence of a filtration.
\end{proof}

\begin{rem}
The above spectral sequence is a special case of a more general Theorem
\cite[ App. II]{g-z}. The above proof
is taken from \cite[Chap. I]{ma1} where the functorial behavior of the
spectral sequence is more visible. For more details see \cite{ma1}.
\end{rem}

\begin{dfn}
Let $X$ be a poset. A map  $\hht_X : X \arr \z_{\ge 0}$ is called
height function if it is a strictly increasing map.
\end{dfn}

For example the height function
$\hht_X(x) = 1 +{\rm dim}(Link_X^-(x))$ is the usual one considered in
\cite{quil7}, \cite{ma1} and \cite{ch}.

\begin{lem}\label{char}
Let X be a poset such that $\lii$ is $(n-{\rm ht}_X(x)-2)$-acyclic,
for every $x \in X$, where $\hht_X$ is a height function on $X$.
Let $\ff: X \arr  \underline{{\rm Ab}}$ be a functor such
that $\ff(x)=0$ for all $x \in X$ with ${\rm ht}_X(x) \ge m$, where
$m \ge 1$. Then
$H_k(X, \ff) = 0$ for $k \le n-m$.
\end{lem}
\begin{proof}
First consider the case of a functor $\ff$ such that $\ff(x)=0$ if
${\rm ht}_X(x) \neq m-1$. Then
$C_k(X, \ff)=\underset{\underset{{\rm ht}_X(x_0)=m-1}{x_0 <\dots
<x_k}}{\bigoplus}
\ff(x_0)$.
Clearly $0=d^k_0=\ff(x_0 < x_1)=\ff(x_0) \arr \ff(x_1)$.
Thus $\partial_k = \s^k_{i=1}(-1)^id^k_i$.
Define $C_{-1}(Link^+_X(x_0), \ff(x_0))=\ff(x_0)$ and complete the singular
complex of $Link^+_X(x_0)$ with coefficient in $\ff(x_0$) to
\[
\cdots \arr C_0(Link^+_X(x_0), \ff(x_0)) \overset{\varepsilon}{\arr}
C_{-1}(Link^+_X(x_0), \ff(x_0)) \arr 0
\]
where $\varepsilon((g_i))=\s_ig_i$. Then
\begin{gather*}
\!\!\!\!\!\!\!\!\!\!\!\!\!\!\!\!\!\!\!\!\!C_k(X, \ff)=
\underset{{\rm ht}_X(x_0)=m-1}{\bigoplus}(\underset
{\underset{x_0<x_1}{x_1 <\dots <x_k}}
{\bigoplus}\ff(x_0))\\
\ \ \ \ \ \ \ \ \ \ \ \ \ \ \ =\underset{{\rm ht}_X(x_0)=m-1}
{\bigoplus}C_{k-1}(Link^+_X(x_0), \ff(x_0)).
\end{gather*}
The complex $C_{k-1}(Link^+_X(x_0), \ff(x_0))$ is the standard
complex for computing the reduced homology of $Link^+_X(x_0)$ with
constant coefficient $\ff(x_0)$.
So
\[
H_k(X,\ff)=\underset{{\rm ht}_X(x)=m-1}{\bigoplus}\H_{k-1}(Link^+_X(x),
\ff(x)).
\]
If ${\rm ht}_X(x_0)=m-1$ then $Link^+_X(x_0)$ is $(n-(m-1)-2)$-acyclic, and by
the universal coefficient theorem \cite[Chap. 5, Thm. 8]{spa},
$\H_{k-1}(Link^+_X(x_0), \ff(x_0))=0$ for $-1 \le k-1 \le n-(m-1)-2$.
This shows that $H_k(X, \ff)=0$ for $0 \le k \le n-m$.
To prove the lemma in general,
we argue by induction on $m$. If $m=1$ then for ${\rm ht}_X(x)\ge 1$,
$\ff(x)=0$. So the lemma follows from the special case above.
Suppose $m \ge 2$. Define $\ff_0$ and  $\ff_1$ to be the functors
\begin{gather*}
\ff_0(x)= \begin{cases}\ff(x) & \text{if ${\rm ht}_X(x)<m-1$}\\
0 & \text{if ${\rm ht}_X(x)\ge m-1$} \end{cases}
,\
\ff_1(x)= \begin{cases}\ff(x) & \text{if ${\rm ht}_X(x)=m-1$}\\
0 & \text{if ${\rm ht}_X(x)\neq m-1$} \end{cases}
\end{gather*}
respectively. Then there is a short exact sequence
$0 \arr \ff_1 \arr \ff \arr \ff_0 \arr 0$.
By the above discussion, $H_k(X, \ff_1)=0$ for
$0\le k \le n-m$ and by induction
for $m-1$, we
have $H_k(X, \ff_0)=0$ for $k \le n-(m-1)$. By the long exact sequence
for the above short exact sequence of functors it is easy to see  that
$H_k(X, \ff)=0$ for $0 \le k \le n-m$.
\end{proof}

\begin{thm}\label{quil}
Let $f :X \arr Y$ be a map of posets and $\hht_Y$ a height function
on $Y$. Assume for every $y \in Y$, that $Link^+_Y(y)$ is
$(n-{\rm ht}_Y(y)-2)$-acyclic and  $f/y$ is $({\rm ht}_Y(y)-1)$-acyclic.
Then $f_\ast: H_k(X, \z) \arr H_k(Y, \z)$ is an
isomorphism for $0 \le k \le n-1$.
\end{thm}
\begin{proof}
By theorem \ref{g-z}, we have the first quadrant spectral sequence
\[
E_{p, q}^2 = H_p(Y, y \mt H_q(f/y, \z)) \Rightarrow H_{p+q}(X, \z).
\]
Since $H_q(f/y, \z)=0$ for $0 < q \le {\rm ht}_Y(y) -1$, the functor
$\gg_q : Y \arr \underline{\rm Ab}$, $\gg_q(y)=H_q(f/y, \z)$ is trivial
for ${\rm ht}_Y(y) \ge q+1$, $ q >0$. By lemma \ref{char},
$H_p(Y, \gg_q)=0$ for $p \le n-(q+1)$. Hence $E^2_{p, q}=0$ for
$p+q \le n-1$, $ q >0$. If $q=0$, by writing the long exact sequence for the
short exact sequence $0 \arr \H_0(f/y, \z) \arr H_0(f/y, \z) \arr \z \arr 0$,
valid because $f/y$ is nonempty, we have
\begin{gather*}
\cdots \arr H_n(Y, \z) \arr H_{n-1}(Y, y \mt \H_0(f/y, \z)) \arr E^2_{n-1, 0}\\
\arr H_{n-1}(Y, \z) \arr H_{n-2}(Y, y \mt \H_0(f/y, \z)) \arr E^2_{n-2,
0}\arr\\
\cdots \arr H_1(Y, \z) \arr H_0(Y, y \mt \H_0(f/y, \z)) \arr E^2_{0, 0} \arr
H_0(Y, \z) \arr 0.
\end{gather*}
If ${\rm ht}_Y(y) \ge 1$, then $ \H_0(f/y, \z)=0$. By lemma \ref{char},
$H_k(Y, y \mt \H_0(f/y, \z))=0$ for $0 \le k \le n-1$.
Thus
\[
E^2_{p, q}= \begin{cases}H_p(Y, \z) & \text{if $q=0$, $0 \le p \le n-1$}\\
0 & \text{if $p+q \le n-1$, $ q>0 $} \end{cases}.
\]
This shows that
$E^2_{p, q}\simeq E^3_{p, q}\simeq \dots \simeq E^r_{p, q}\simeq
\dots \simeq E^\infty_{p, q}$ for
$0 \le p+q \le n-1$.
For every $k$, we have a
filtration
$0=F_{-1}H_k \se F_0 H_k \se \dots \se F_kH_k=H_k(X, \z)$ of $H_k(X, \z)$,
such that $ E^\infty_{p, q}\simeq  F_p H_{p+q}/F_{p-1}H_{p+q}$
\cite[Chap. 5, 5.2.6]{weib}. Let $0 \le k \le n-1$.
For $0 \le i < k$,  we have $0=  E^\infty_{i, k-i}\simeq  F_iH_k/F_{i-1}H_k$,
so $F_iH_k=F_{i-1}H_k$. Hence
$0=F_{-1}H_k = F_0 H_k= \dots = F_{k-1}H_k$.
On the other hand $E^\infty_{k, 0}\simeq  F_kH_k/F_{k-1}H_k$.
Therefore $H_k(X, \z) \simeq H_k(Y, \z)$. This shows that
$H_k(X, \z) \simeq H_k(Y, \z)$ for $0 \le k \le n-1$.
Now consider the commutative diagram
\begin{gather*}
\begin{array}{ccc}
X  & \overset{f}{\longrightarrow} & Y    \\
\Big\downarrow\vcenter{%
\rlap{$\scriptstyle{f}$}} &      &
\Big\downarrow\vcenter{%
\rlap{$\scriptstyle{{\rm id}_Y}$}}\\
Y          & \overset{{\rm id}_Y}{\longrightarrow}  &  Y .
\end{array}
\end{gather*}
By functoriality of the spectral sequence, and the above calculation
we get the diagram
\begin{gather*}
\begin{array}{ccc}
H_k(Y, y \mt H_0(f/y,\z))& \overset{\simeq}{\longrightarrow} & H_k(X, \z)    \\
\Big\downarrow\vcenter{%
\rlap{$\scriptstyle{{{\rm id}_Y}_\ast}$}} &      &
\Big\downarrow\vcenter{%
\rlap{$\scriptstyle{f_\ast}$}}\\
H_k(Y, y \mt H_0({\rm id}_Y/y,\z))&
\overset{\simeq}{\longrightarrow}  & H_k(Y, \z).
\end{array}
\end{gather*}
Since ${\rm id}_Y/y=Y_{\le y}$ is contractible, we have
$H_k(Y, y \mt H_0({\rm id}_Y/y,\z)) = H_k(Y, \z)$. The map ${{\rm id}_Y}_\ast$
is an isomorphism for $0 \le k \le n-1$, from the
above long exact sequence. This shows that $f_\ast$ is an isomorphism for
$0 \le k \le n-1$.
\end{proof}

\begin{lem}\label{quil2}
Let X be a $0$-connected poset. Then X is $1$-connected if and only if
for every local system $\ff$ on X and every $x \in X$, the map
$ \ff(x) \arr H_0(X, \ff)$, induced from the inclusion $\{x\} \harr X$,
is an isomorphism $($or equivalently, every local system on X is
a isomorphic with a constant local system$)$.
\end{lem}
\begin{proof}
If $X$ is 1-connected then by theorem \ref{wh1} and the connectedness of $X$,
one has
$\ff(x) \overset{\simeq}{\longrightarrow}H_0(X, \ff)$
for every $x \in X$. Now let every local system on $X$ be isomorphic
with a constant local system. Let $\ff: X \arr \underline{{\rm Set}}$
be in $\lll_S(X)$. Define the functor
$\gg: X \arr \underline{{\rm Ab}}$ where $\gg(x)$
is the free abelian group generated by $\ff(x)$. Clearly $\gg$ is a
local system  and so it is constant system. This shows that $\ff$
is isomorphic to a constant functor. So by lemma \ref{g-z-q}, any
connected  covering space of $|X|$ is isomorphic to $|X|$. This shows that
the universal covering of $|X|$, is $|X|$. Note that the
universal covering of a connected simplicial simplex exists
and is simply connected \cite[Chap. 2, Cor. 14 and 15]{spa}.
Therefore $X$ is 1-connected.
\end{proof}

\begin{thm}\label{quil4}
Let $f :X \arr Y$ be a map of  posets and  $\hht_Y$ a height function
on $Y$. Assume for every $y \in Y$, that $Link^+_Y(y)$
is $(n-{\rm ht}_Y(y)-2)$-connected and $f/y$ is $({\rm ht}_Y(y)-1)$-connected.
Then $X$ is $(n-1)$-connected if and only if $Y$ is $(n-1)$-connected.
\end{thm}
\begin{proof}
By \ref{hur} and \ref{quil} it is enough to prove that $X$ is 1-connected
if and only if $Y$ is 1-connected, when $n \ge 2$.
Let $\ff : X \arr \underline{\rm Ab}$
be a local system. Define the functor $\gg : Y \arr  \underline{\rm Ab}$
with
\[
\gg(y)=\begin{cases}H_0(f/y, \ff) & \text{if ${\rm ht}_Y(y) \neq 0$}\\
H_0(Link^+_Y(y), y' \mt H_0(f/y', \ff)) & \text{if ${\rm ht}_Y(y)=0$}
\end{cases}.
\]
We prove that $\gg$ is a local system. If ${\rm ht}_Y(y) \ge 2$ then $f/y$ is
1-connected and by \ref{quil}, $\ff|_{f/y}$ is a constant system, so
by \ref{quil2}, $H_0(f/y, \ff) \simeq \ff(x)$ for every $x\in f/y$. If
${\rm ht}_Y(y)=1$, then $f/y$ is 0-connected and $ Link^+_Y(y)$ is nonempty.
Choose $y' \in Y$ such that $y<y'$. Now $f/y'$ is 1-connected and so
$\ff|_{f/y'}$ is a constant system on $f/y'$. But $f/y \subset f/y'$,
so $\ff|_{f/y}$ is a constant
system. Since $f/y$ is 0-connected, by \ref{wh1} and the fact that we
mentioned before theorem \ref{q-w},
$H_0(f/y, \ff) \simeq \ff(x)$ for every $x \in f/y$. Now let
${\rm ht}_Y(y)=0$. Then $Link^+_Y(y)$ is 0-connected, $f/y$
is nonempty and for every
$y' \in Link^+_Y(y)$, $H_0(f/y', \ff) \simeq H_0((f/y)^\circ, \ff)$
where $(f/y)^\circ$ is  a component of $f/y$, which we fix.
This shows that the
local system  $\ff': Link^+_Y(y) \arr \underline{\rm Ab}$ with
$y' \mt H_0(f/y', \ff)$ is isomorphic to a constant system, so
$H_0(Link^+_Y(y), y' \mt H_0(f/y', \ff))=H_0(Link^+_Y(y), \ff')
\simeq \ff'(y') \simeq \ff(x)$ for every $x \in f/y'$. This shows that
$\gg$ is a local system.

If $Y$ is 1-connected, by \ref{quil2}, $\gg$ is a constant system.
But it is easy to see that $\ff \simeq \gg \circ f$. Therefore $\ff$
is a constant system. Since $X$ is connected by our homology calculation,
by \ref{quil2} we conclude that $X$ is 1-connected.
Now let $X$ be 1-connected. If $\ee$ is a local system on $Y$, then
$f^\ast\ee:= \ee \circ f$ is a local system on $X$. So it is a
constant local system. As above we can construct a local
system $\gg'$ on $Y$ from $\ff':=\ee \circ f$. This gives a
natural transformation from $\gg'$ to $\ee$ which is an isomorphism.
Since $\ee \circ f$ is constant, by \ref{wh1} and \ref{quil2}
and an argument as above one sees that $\gg'$ is constant.
Therefore $\ee$ is isomorphic to
a constant local system and \ref{quil2} shows that $Y$ is
$1$-connected.
\end{proof}


\begin{rem}\label{equiv}
In the proof of the above theorem \ref{quil4} we showed in fact that:
Let $f :X \arr Y$ be a map of  posets and  $\hht_Y$ a height function
on $Y$. Assume for every $y \in Y$, that $Link^+_Y(y)$
is $(n-{\rm ht}_Y(y)-2)$-connected and $f/y$ is $({\rm ht}_Y(y)-1)$-connected.
Then $f^\ast: \lll_S(Y) \arr \lll_S(X)$, with $\ee \mt \ee \circ f$ is an
equivalence of categories.
\end{rem}

\begin{rem}
Theorem \ref{quil4} is a generalization of a theorem of Quillen
\cite[Thm. 9.1]{quil7}.
We proved that the converse of that theorem is also valid. Our proof is
similar in outline
to the  proof by Quillen. Furthermore,
lemma \ref{char} is a generalized
version of lemma 1.3 from  \cite{ch}. With more restrictions, Maazen,
in \cite[ Chap. II]{ma1} gave an easier proof of  Quillen's theorem.
\end{rem}

\section{Homology and homotopy of posets of sequences}

Let $V$ be a nonempty set. We denote by $\o(V)$ the poset of finite
ordered sequences of
distinct elements of $V$, the length of each sequence being at least one.
The partial ordering on $\o(V)$ is defined by refinement: $(v_1, \dots ,v_m)
\le (w_1, \dots , w_n)$ if and only if there is a strictly increasing map
$\phi : \{ 1, \dots , m\} \arr \{ 1, \dots , n\}$ such that
$v_i = w_{\phi(i)}$, in other words, if $(v_1, \dots ,v_m)$ is an
order preserving subsequence of $(w_1, \dots , w_n)$.
If $v = (v_1, \dots ,v_m)$
we denote by $|v|$ the length of $v$, that is $|v|=m$. If
$v = (v_1, \dots ,v_m)$ and $w = (w_1, \dots ,w_n)$, we write
$(v_1, \dots ,v_m, w_1, \dots ,w_n)$ as $vw$. Define $F_v$ to be the set
of $w \in F$  such that $wv \in F$. Note that ${(F_v)}_w = F_{wv}$.
A subset $F$ of $\o(V)$ is said to satisfy the {\it chain condition}
if  $v \in F$ whenever $w \in F$, $v \in \o(V)$ and $v \le w$.
The subposets of $\o(V)$ which satisfy the chain condition are extensively
studied in \cite{ma1}, \cite{kall} and \cite{ch2}. In this section we will
study them some more.

Let $F \se \o(V)$. For a nonempty set $S$ we define the poset $\fs$ as
\[
\fs :=\{ ((v_1, s_1),\dots, (v_r, s_r)) \in \o(V \times S) :
(v_1,\dots, v_r) \in F\}.
\]
Assume $s_0 \in S$ and consider the injective poset map
$l_{s_0}: F \arr \fs$ with
$(v_1,\dots, v_r) \mt ((v_1, s_0),\dots, (v_r, s_0))$.
We have clearly a projection $p: \fs \arr F$ with
$((v_1, s_1),\dots, (v_r, s_r)) \mt (v_1,\dots, v_r)$ such
that $p\circ l_{s_0}={\rm id}_F$.

\begin{lem}\label{maazen5}
Suppose $F \se \o(V)$ satisfies the chain condition and $S$ is a
nonempty set. Assume for every $v \in F$, that $F_v$ is $(n-|v|)$-connected.
\par {\rm (i)} If $s_0 \in S$ then
$(l_{s_0})_\ast: H_k(F, \z) \arr H_k(\fs, \z)$ is an isomorphism for
$0 \le k \le n$.
\par {\rm (ii)} If F is ${\rm min}\{1, n-1\}$-connected, then
 $(l_{s_0})_\ast: \pi_k(F, v) \arr \pi_k(\fs, l_{s_0}(v))$ is an
isomorphism for $0 \le k \le n$.
\end{lem}
\begin{proof}
This follows by \cite[Prop. 1.6]{ch2} from the fact that
$p \circ l_{s_0}={\rm id}_F$.
\end{proof}

\begin{lem}\label{maazen1}
Let $F \se \o(V)$  satisfies the chain condition.
Then $|Link_F^-(v)| \simeq S^{|v|-2}$ for every $v \in F$.
\end{lem}
\begin{proof}
Let $v=(v_1,\dots ,v_n)$. By definition
$Link_F^-(v) =\{w \in F: w < v\} =
\{(v_{i_1},\dots ,v_{i_k}) : k<n, i_1< \dots < i_k\}$. Hence $|Link_F^-(v)|$
is isomorphic to the barycentric subdivision of the boundary
of the standard simplex
$\dl_{n-1}$. It is well
known that $\partial\dl_{n-1} \simeq S^{n-2}$, hence $|Link_F^-(v)| \simeq
S^{|v|-2}$.
\end{proof}

\begin{thm}[Nerve Theorem for Posets]\label{p-n-t}
Let $V$ and $T$ be two nonempty sets, $F \se \o(V)$ and $X \se \o(T)$.
Assume $X=\bcu_{v \in F}X_v$ such that if $v \le w$ in $F$, then $X_w \se X_v$.
Let $F$, $X$ and $X_v$, for every $v \in F$, satisfy the chain condition.
Also assume
\par {\rm (i)} for every $v \in F$, $X_v$ is $(l-|v|+1)$-acyclic
$($resp. $(l-|v|+1)$-connected$)$,
\par {\rm (ii)} for every $x \in X$, $\a_x:=\{ v \in F: x \in X_v\}$
is $(l-|x|+1)$-acyclic $($resp. $(l-|x|+1)$-connected$)$.\\
Then $H_k(F, \z) \simeq H_k(X, \z)$ for $0 \le k \le l$
$($resp. $F$ is l-connected if and only if X is l-connected$)$.
\end{thm}
\begin{proof}
Let $F_{\le l+2}=\{ v \in F : |v| \le l+2\}$ and let $i: F_{\le l+2} \arr F$
be the inclusion. Clearly $|F_{\le l+2}|$ is the $(l+1)$-skeleton of
$|F|$, if we consider $|F|$ as a cell complex whose $k$-cells are
the $|F_{\leq v}|$ with $|v|=k+1$. It is well known that
$i_\ast: H_k(F_{\le l+2}, \z) \arr H_k(F, \z)$ and
$i_\ast: \pi_k(F_{\le l+2}, v) \arr \pi_k(F, v)$ are isomorphisms for
$0 \le k \le l$ (see \cite{whi}, Chap. II, corollary 2.14, and
\cite{whi}, Chap. II, Corollary 3.10 and Chap. IV lemma 7.12.)
So it is enough to prove the theorem for $F_{\le l+2}$ and $X_{\le l+2}$.
Thus assume $F=F_{\le l+2}$ and $X=X_{\le l+2}$. We define
$Z \se X \times F$ as $Z=\{(x, v): x \in X_v\}$. Consider the projections
\begin{gather*}
f: Z \arr F,\; (x, v) \mt v \qquad , \qquad g: Z \arr X, \;(x, v) \mt x.
\end{gather*}
First we prove that $f^{-1}(v) \sim v \bs f$ and $g^{-1}(x) \sim x \bs g$,
where $\sim$ means homotopy equivalence.
By definition $v \bs f=\{(x, w): w \ge v,\; x \in X_w\}$.
Define $\phi: v \bs f \arr f^{-1}(v)$, $(x, w) \mt (x, v)$.
Consider the inclusion $j: f^{-1}(v) \arr v \bs f$. Clearly
$\phi\circ j(x, v)=\phi(x, v)=(x, v)$ and
$j \circ \phi(x, w)=j(x, v)=(x, v) \le (x, w)$. So by
\ref{qui}(ii), $v \bs f$ and $f^{-1}(v)$ are homotopy equivalent.
Similarly $ x \bs g \sim g^{-1}(x)$.

Now we prove that the maps $f^{op} : Z^{op} \arr Y^{op}$ and
$g^{op} : Z^{op} \arr X^{op}$ satisfy the conditions of \ref{quil}.
First  $f^{op} : Z^{op} \arr Y^{op}$; define the height
function  $\hht_{F^{op}}$ on $\fop$ as $\hht_{F^{op}}(v)=l+2-|v|$.
It is easy to see that
$f^{op}/v \simeq v \bs f \sim f^{-1}(v) \simeq X_v$.
Hence $f^{op}/v$ is $(l-|v| +1)$-acyclic (resp. $(l-|v| +1)$-connected).
But $l-|v| +1=(l+2-|v|)-1=\hht_{F^{op}}(v)-1$, so $f^{op}/v$ is
$(\hht_{F^{op}}(v)-1)$-acyclic (resp. $(\hht_{F^{op}}(v)-1)$-connected).
Let $n:=l+1$. Clearly $Link^+_{F^{op}}(v)=Link^-_F(v)$. By lemma
\ref{maazen1}, $|Link^-_F(v)|$ is $(|v|-3)$-connected.
But  $(|v|-3)=l+1-(l+2-|v|)-2=n-\hht_{F^{op}}(v)-2$. Thus $Link^+_{F^{op}}(v)$
is $(n-\hht_{F^{op}}(v)-2)$-acyclic (resp. $(n-\hht_{F^{op}}(v)-2)$-connected).
Therefore by theorem \ref{quil}, $f_\ast: H_i(Z, \z) \arr H_i(F, \z)$ is an
isomorphism for $0 \le i \le l$ (resp. by \ref{quil4}, $F$ is $l$-connected
if and onl if $Z$ is $l$-connected).
Now consider $g^{op} : Z^{op} \arr X^{op}$. We saw in the above that
$g^{op}/x \simeq x \bs g \sim g^{-1}(x)$ and
$g^{-1}(x)=\{(x, v): x \in X_v\}
\simeq \{ v \in F:  x \in X_v\}$. It is  similar to the
case of $f^\op$ to see that $g^\op$
satisfies the conditions of theorem \ref{quil}, hence
$g_\ast: H_i(Z, \z) \arr H_i(X, \z)$ is an
isomorphism for $0 \le i \le l$ (resp. by \ref{quil4}, $X$ is $l$-connected
if and only if $Z$ is $l$-connected). This completes the proof.
\end{proof}

Let $K$ be a simplicial complex and $\{K_i\}_{i \in I}$ a family of
subcomplexes such that $K=\bigcup_{i \in I}K_i$. The {\it nerve} of this
family of subcomplexes of $K$  is the simplicial complex $\n(K)$
on the vertex set $I$ so that a finite subset
$\si \se I$ is in $\n(K)$ if and only if
$\bigcap_{i \in \si}K_i \neq \varnothing$. The nerve $\n(K)$ of $K$,
with the inclusion relation,
is a poset. As we already said we can consider a simplicial complex
as a poset of its simplices. We leave it
to the interested reader to derive the next result
from the above theorem.

\begin{cor}[Nerve Theorem]\label{h-n}
Let $K$ be a simplicial complex and $\{ K_i \}_{i \in I}$ a family of
subcomplexes such that $K=\bigcup_{i \in I}K_i$. Suppose every nonempty
finite intersection $\bigcap_{j=1}^tK_{i_j}$ is $(l-t+1)$-acyclic
$($resp. $(l-t+1)$connected$)$.
Then $H_k(K, \z) \simeq H_k(\n(K), \z)$ for $0 \le k \le l$
$($resp. $K$ is $l$-connected if and only if  $\n(K)$ is $l$-connected$)$.
\end{cor}

\begin{rem}
In \cite{gra}, a special case of the theorem \ref{p-n-t} is proved.
The nerve theorem for a
simplicial complex \ref{h-n}, in the stated generality, is
proved for the first time in \cite{bjo1}, see also \cite[ p.~1850]{bjo}.
For more information about different types of nerve theorem and
more references about them see \cite[ p.~1850]{bjo}.
\end{rem}

\begin{lem}\label{h0}
Let $F \se \o(V)$ satisfy the chain condition and let
$\gg: F^{op} \arr \underline{{\rm Ab}}$ be a functor. Then the natural
map $\psi:\bigoplus_{v \in F,\; |v|=1}\gg(v) \arr H_0(F^{op}, \gg)$ is
surjective.
\end{lem}
\begin{proof}
By definition $C_0(\fop, \gg)=\bigoplus_{v \in \fop}\gg(v)$,
$C_1(\fop, \gg)=\bigoplus_{v<v' \in \fop}\gg(v)$ and we have the chain complex
\begin{gather*}
\cdots \arr C_1(\fop, \gg) \overset{\partial_1}{\arr} C_0(\fop, \gg) \arr 0,
\end{gather*}
where $\partial_1=d_0^1 - d_1^1$. Again by definition
$H_0(F^{op}, \gg)=C_0(\fop, \gg)/\partial_1$. Now let $w \in F$ and $|w|\ge 2$.
Then there is a $v \in F$, $v\leq w$,
with $|v|=1$. So $w < v$ in $\fop$, and
we have the component
$\partial_1|_{\gg(w)}: \gg(w) \arr \gg(w)\oplus \gg(v)$,
$x \mt d_0^1(x)-d_1^1(x)=d_0^1(x)-x$. This shows that
$\gg(w) \se {\rm im}\partial_1+{\rm im}\psi$. Therefore $ H_0(F^{op}, \gg)$
is generated by the groups $\gg(v)$ with $|v|=1$.
\end{proof}

\begin{thm}\label{surj}
Let $V$ and $T$ be two nonempty sets, $F \se \o(V)$ and $X \se \o(T)$.
Assume $X=\bcu_{v \in F}X_v$ such that if $v \le w$ in $F$, then $X_w \se X_v$
and let $F$, $X$ and $X_v$, for every $v \in F$, satisfy the chain condition.
Also assume
\par {\rm (i)} for every $v \in F$, $X_v$ is
${\rm min}\{l-1, l-|v|+1\}$-connected,
\par {\rm (ii)} for every $x \in X$, $\a_x:=\{ v \in F: x \in X_v\}$
is $(l-|x|+1)$-connected,
\par {\rm (iii)} $F$ is $l$-connected.\\
Then $X$ is $(l-1)$-connected and the natural map
\begin{gather*}
\bigoplus_{v \in F,\;|v|=1}(i_v)_\ast: \bigoplus_{v \in F,\;|v|=1}H_l(X_v, \z)
\arr H_l(X, \z)
\end{gather*}
is surjective, where $i_v:X_v \arr X$ is the inclusion.
Moreover, if for every $v$ with $|v|=1$, there is
an $l$-connected $Y_v$ with $X_v \se Y_v\se X$, then $X$ is also $l$-connected.
\end{thm}
\begin{proof}
If $l=-1$, then everything is easy. If $l=0$, then for $v$ of length
one, $X_v$ is nonempty, so $X$ is nonempty. This shows that $X$ is
$(-1)$-connected. Also,
every connected component of
$X$ intersects at least one $X_w$ and therefore also contains
a connected component of an $X_v$ with $|v|=1$.
This gives the surjectivity
of the homomorphism
\begin{gather*}
\bigoplus_{v \in F,\; |v|=1}(i_v)_\ast: \bigoplus_{v \in F,\; |v|=1}H_0(X_v,
\z)
\arr H_0(X, \z).
\end{gather*}
Now assume that, for every $v$ of length one, $X_v \se Y_v$
where $Y_v$ is connected. We prove, in a combinatorial way, that
$X$ is connected. Let $x, y \in X$, $x \in X_{(v_1)}$ and
$y \in X_{(v_2)}$ where $(v_1), (v_2) \in F$. Since $F$ is connected,
there is a sequence $(w_1), \dots , (w_r) \in F$ such that they give a path,
in $F$, from $(v_1)$ to $(v_2)$, that is
\begin{gather*}
\begin{array}{ccccccccccc}
 (v_1) &          &         &       &\!\!\!\!\!\!(w_1)&     & (w_r) &
&           &         &\!\!\!\!\!(v_2)  \\
      &\!\!\!\!\!\!\diagdown &         &\!\!\!\diagup&     &\dots&
& \!\!\!\diagdown &           & \!\!\!\!\!\!\diagup &       \\
      &          &\!\!\!(v_1,w_1)&       &     &     &       &
& \!\!\!(w_r, v_2)&         & .
\end{array}
\end{gather*}
Since $Y_{(v_1)}$ is connected, $x \in X_{(v_1)} \se Y_{(v_1)}$ and
$X_{(v_1, w_1)} \neq \varnothing$, there is an element $x_1 \in X_{(v_1, w_1)}$
such that there is a path, in $Y_{(v_1)}$, from $x$ to $x_1$.
Now $x_1 \in Y_{(w_1)}$.
Similarly we can find $x_2 \in X_{(w_1, w_2)}$ such that
there is a path, in $Y_{(w_1)}$,
from $x_1$ to $x_2$. Now $x_2 \in Y_{(w_2)}$.
Repeating this process finitely many times, we find a path from $x$ to
$y$. So $X$ is connected.

Hence we assume that $l \ge 1$. As we said in the proof of
theorem \ref{p-n-t}, we can assume that
$F= F_{\le l+2}$ and $X=X_{\le l+2}$ and we define $Z$, $f$ and $g$ as
we defined them there. Define the height function $\hht_{F^{op}}$ on
$\fop$ as $\hht_{F^{op}}(v)=l+2-|v|$.
As we proved in the proof of theorem \ref{p-n-t},
$f^{op}/v \simeq v\bs f \sim f^{-1}(v) \simeq X_v$. Thus $f^{op}/v$ is
$(\hht_{F^{op}}(v)-1)$-connected if $|v|>1$ and it is
$(\htf-2)$-connected if $|v|=1$ and also
$|Link_{\fop}^+(v)|$ is $(l+1-\hht_{F^{op}}(v)-2)$-connected.
By theorem \ref{g-z}, we have the first quadrant spectral sequence
\begin{gather*}
E_{p, q}^2 = H_p(\fop, v \mt H_q(f^\op/v, \z)) \Rightarrow H_{p+q}(Z^\op, \z).
\end{gather*}
For $0 <q \le \htf-2$, $H_q(f^\op/v, \z)=0$. Define
$\gg_q : \fop \arr \underline{\rm Ab}$, $\gg_q(v)=H_q(f^\op/v, \z)$.
Then $\gg_q(v)=0$ for $\htf \ge q+2$, $ q >0$. By lemma \ref{char},
$H_p(\fop, \gg_q)=0$ for $p \le l+1-(q+2)$. Therefore $E^2_{p, q}=0$ for
$p+q \le l-1$, $ q >0$. If $q=0$, arguing similarly to the
proof of theorem \ref{quil},
we get $E_{p, 0}^2=0$ if $0<p\le l-1$ and $E_{0, 0}^2= \z$. Also by the
fact that $\fop$ is $l$-connected we get the surjective homomorphism
$ H_l(\fop, v \mt \H_0(f^\op/v, \z)) \two E^2_{l, 0}$.
Since $l \ge 1$, $\H_0(f^\op/v, \z)=0$ for all $v \in \fop$
with $\hht_{F^{op}}(v)\geq1$ and so
$H_l(\fop, v \mt \H_0(f^\op/v, \z))=0$
by lemma \ref{char}. Therefore $E^2_{l, 0}=0$.
Let ${\gg'}_q : \fop \arr \underline{\rm Ab}$,
${\gg'}_q(v)=\begin{cases}0 & \text{if $\htf<l+1$}\\
H_q(f^\op/v, \z) & \text{if $\htf=l+1$} \end{cases}$  and
${\gg''}_q : \fop \arr \underline{\rm Ab}$,
${\gg''}_q(v)= \begin{cases}H_q(f^\op/v, \z) & \text{if $\htf<l+1$}\\
0 & \text{if $\htf=l+1$} \end{cases}$. Then we have the
short exact sequence
$0 \arr {\gg'}_q \arr {\gg}_q \arr {\gg''}_q \arr 0$
and the associated long exact sequence
\begin{gather*}
\!\!\!\!\!\!\!\!\!\cdots \arr H_{l-q}(\fop, {\gg'}_q)
\arr H_{l-q}(\fop, \gg_q) \arr \\
\ \ \ \ \ \ \ \ \ \ \ \ \ H_{l-q}(\fop, {\gg''}_q) \arr
H_{l-q-1}(\fop, {\gg'}_q) \arr \cdots.
\end{gather*}
If $q>0$, then ${\gg''}_q(v)=0$ for $0 < q \le \htf-1$ and so by lemma
\ref{char}, $H_p(\fop, {\gg''}_q)=0$ for $p+q \le l$, $q>0$. Also if $|v|=1$
then $H_0(f^\op/v, \z)=0$ for $0 <q \le \htf-2=l-1$.
This shows ${\gg'}_q=0$ for $0 < q \le l-1$. From the long
exact sequence  and the above calculation we get,
$E_{p, q}^2= \begin{cases}\z & \text{if $p=q=0$}\\
0 & \text{if $0<p+q \le l$, $q \neq l$} \end{cases}$.
\begin{gather*}
\begin{array}{ccccccl}
\setlength{\unitlength}{1em}
\begin{picture}(0,0)(-10.75,8.3)
\put(-4.65,6.2){\vector(-2,1){2.5}}
\put(-8.5,1){\line(0,1){8}}
\put(-8.5,1){\line(1,0){12}}
\end{picture} l+1 &  * & & & & & E^2_{p, q}\\
l & \ast  &      * & & & & \\
 &  0 &       0   &   * &  & & \\
 & \vdots & \vdots  & \ddots & \ddots &  &
\\  &  0     &  \cdot  & \cdot  &    \ddots & * & \\
0 & \mathbb Z  &  0  &    \cdots   &   0    &     0& *\\
& 0 & 1 &  &  &  \ l & l+1
\end{array}
\end{gather*}
Thus for $0 \le p+q \le l$, $q \neq l$,
$E_{p, q}^2 \simeq E_{p, q}^3\simeq \dots \simeq E_{p, q}^\infty$ and
we have a filtration of $H_l(Z^\op, \z)$,
$0=F_{-1}H_l \se F_0 H_l \se \dots \se F_lH_l=H_l(Z^\op, \z)$
such that $E^\infty_{p, q}\simeq  F_pH_{p+q}/F_{p-1}H_{p+q}$.
If $i \neq 0$ then $0=  E^\infty_{i, l-i}\simeq  F_iH_l/F_{i-1}H_l$,
so $F_iH_k=F_{i-1}H_k$. Therefore
$0=F_{-1}H_l \se F_0 H_l= \dots = F_lH_l=H_l(Z^\op, \z)$ and
$E^\infty_{0, l}\simeq  F_0H_l/F_{-1}H_l \simeq H_l(Z^\op, \z)$.
By definition
$E_{0, l}^{r+1}={\rm ker}(d_{0, l}^r)/{\rm im}(d_{r, l-r+1}^r)$. Thus there
exist $r \in \z$ such that
$E_{0, l}^2 \two E_{0, l}^3\two \dots \two E_{0, l}^r \simeq E_{0, l}^{r+1}
\simeq \dots \simeq E_{0, l}^\infty.$
Hence we get a surjective map
$H_0(\fop, v \mt H_l(f^\op/v, \z)) \two H_l(Z^\op, \z).$
By lemma \ref{h0}, we have a surjective map
$\bigoplus_{v \in F,\; |v|=1}H_l(f^\op/v, \z) \two H_l(Z^\op, \z).$

Now consider the map $g^\op :Z^\op \arr X^\op$ and define the height function
$\hht_{X^\op}(x)=l+2-|x|$ on $X^\op$. Arguing similarly to the proof of
theorem \ref{p-n-t} one sees that
$g_\ast: H_k(Z, \z) \arr H_k(X, \z)$ is an isomorphism for
$0 \le k \le l$. Therefore we get a surjective map
$\bigoplus_{v \in F,\; |v|=1}H_l(X_v, \z) \two H_l(X, \z)$. We call it $\psi$.
We prove that this map is the same map that we claimed. For $v$ of
length one consider the commutative diagram of posets
\begin{gather*}
\begin{array}{cclccc}
& &\llap{\ $\{v\}^\op$}&
\llap{
\rlap{$\relbar\joinrel\!\relbar\joinrel\!\relbar\joinrel\!
\longrightarrow$}\ \ \ \ }
& &\llap{$\ \ \ \ \ \ \ \ \ \ \ \ \ \ \ \ \fop$}\cr
&\nearrow& & &\!\!\!\!\!\!\!\!\nearrow& \cr
\rlap{$f^{-1}(\{v\})^\op$}\ \ \ & \ \ \ \ \ \
\rlap{\ \ \ $\relbar\joinrel\!\longrightarrow$}
& &\ \ \ \ \ \ \ \ \rlap{\!\!\!\!\!\!$Z^\op$}& &\cr
&\searrow& & &\!\!\searrow& \cr
& &\llap{\ $X_v^\op$}&
\llap{
\rlap{$\relbar\joinrel\!\relbar\joinrel\!\relbar\joinrel\!
\longrightarrow$}\ \ \ \ }
& &\llap{$X^\op$}
\end{array}
\end{gather*}
By functoriality of the spectral sequence for the above diagram and
the lemma \ref{h0} we get the commutative diagram
\begin{gather*}
\begin{array}{ccc}
H_l(f_v^\op/v, \z)& \overset{(j_v)_\ast}{
\relbar\joinrel\!\relbar\joinrel\!\relbar\joinrel\!\relbar\joinrel\!
\longrightarrow}
& \bigoplus_{v\in F,|v|=1}H_l(f^\op/v,\z)\\
 \big\downarrow                   &    &  \big\downarrow \\
H_0(\{v\}^\op, v \mt H_l(f_v^\op/v, \z))  &
\relbar\joinrel\!\relbar\joinrel\!\longrightarrow& H_0(\fop, v \mt
H_l(f^\op/v,\z))\\
     \big\downarrow                     &    &    \big\downarrow\\
 H_l(f^{-1}(v)^\op , \z)  &
\relbar\joinrel\!\relbar\joinrel\!\relbar\joinrel\!\relbar\joinrel\!
\longrightarrow
& H_l(Z^\op, \z)    \\
    \big\downarrow                      &    &  \big\downarrow\\
H_l(X_v^\op, \z)&\overset{(i_v^\op)_\ast}{
\relbar\joinrel\!\relbar\joinrel\!\relbar\joinrel\!\relbar\joinrel\!
\longrightarrow}
&  H_l(X^\op, \z)
\end{array}
\end{gather*}
where $j_v: f_v^\op/v \arr f^\op/v$ is the inclusion which is a
homotopy equivalence as we already mentioned.
It is not difficult to see that the composition of homomorphisms in
the left column of the above diagram
induces the identity map from $H_l(X_v, \z)$,
the composition of homomorphisms in the right
column of above diagram induces the surjective map $\psi$
and the last row induces the homomorphism $(i_v)_\ast$. This show that
$(i_v)_\ast = \psi|_{H_l(X_v, \z)}$. This completes the proof
of surjectiveness.

Now let for $v$ of length one $X_v \se Y_v$ where $Y_v$ is $l$-connected.
Then we have the commutative diagram
\begin{gather*}
\begin{array}{ccc}
H_l(X_v, \z)  &\!\!\!\!\!\!\!\!\!
\overset{(i_v)_\ast}{\relbar\joinrel\!\relbar\joinrel\!\longrightarrow}  &
\!\!\!\!\!\!\!\!\! H_l(X, \z) \\
\ \ \ \ \ \  \searrow &     & \!\!\!\!\!\!\!\!\!\!\!\!\!\!\!\!\!\!\nearrow \\
& \!\!\!\!\!\!\!\!\!\!\!\
H_l(Y_v, \z) &
\end{array}.
\end{gather*}
By the assumption $H_l(Y_v, \z)$ is trivial and
this shows that $(i_v)_\ast$ is the zero map.
Hence by the surjectivity, $H_l(X, \z)$ is trivial. If $l \ge 2$, the nerve
theorem \ref{p-n-t} says that $X$ is simply connected and
by the Hurewicz theorem \ref{hur},
$X$ is $l$-connected. So the only
case that is left is when $l=1$.
By theorem \ref{quil4},
$X$ is 1-connected if and only if $Z$ is 1-connected. So it is enough
to prove that $Z^\op$ is 1-connected. Note that as we said,
we can assume that $F=F_{\le 3}$ and $X=X_{\le 3}$.
Suppose $\ff$ is a local system on $Z^\op$. Define the functor
$\gg : \fop \arr  \underline{\rm Ab}$, as
\[
\gg(y)=\begin{cases}H_0(f^\op/v, \ff) & \text{if $|v|=1, 2$}\\
H_0(Link_{\fop}^+(v), v' \mt H_0(f^\op/ v', \ff)) & \text{if $|v|=3$}
\end{cases}.
\]
We prove that $\gg$ is a local system on $\fop$.
Put $Z_w := g^{-1}(Y_w)$ for $|w|=1$.
If $|v|=1, 2$,
then $f^\op/v$ is 0-connected and $f^\op/v \se Z_w^\op$, where
$w \le v$, $|w|=1$.
By remark \ref{equiv} we can assume that $\ff=\ee \circ g^\op$ where
$\ee$ is a local system on $X^\op$.
Then
$ \ff|_{Z_w^\op}=\ee|_{Y_w^\op} \circ g^\op|_{Z_w^\op}$. Since $Y_w^\op$
is 1-connected, $\ee|_{Y_w^\op}$ is a constant local system. This shows that
$\ff|_{Z_w^\op}$
is a constant local system. So $\ff|_{f^\op/v}$ is a constant local
system and since $f^\op/v$ is 0-connected we have
$H_0(f^\op/v, \z)\simeq \ff(x)$, for every $x \in f^\op/v$.
If $|v|=3$, with an argument similar to the proof of the theorem \ref{quil4}
and the
above discussion one can get
$\gg(v)\simeq \ff(x)$ for every $x \in f^\op/v$. This shows that $\gg$ is a
local system on $\fop$. Hence it is a constant local system, because
$\fop$ is 1-connected.
It is easy to see that $\ff \simeq \gg \circ f$. Therefore $\ff$
is a constant system. Since $X$ is connected by our homology calculation,
by \ref{quil2} we conclude that $X$ is 1-connected. This completes
the proof.
\end{proof}

\section{Posets of isotropic and hyperbolic unimodular sequences}

Let $R$ be an associative ring with unit. A vector
$(r_1, \dots, r_n) \in R^n$ is called unimodular if there exist
$s_1, \dots, s_n \in R$ such that $\s_{i=1}^ns_ir_i=1$, or
equivalently if the submodule generated by this vector is a free
summand of
the right $R$-module $R^n$. We denote the standard basis of
$R^n$ by $e_1, \dots, e_n$. If $n \le m$, we assume that
$R^n$ is the submodule of $R^m$ generated by $e_1, \dots, e_n \in R^m$.

We say that a ring $R$ satisfies the {\it stable range condition}
$({\rm S}_m)$,  if $m\geq1$ is an integer so that
for every unimodular vector $(r_1, \dots ,r_m, r_{m+1}) \in R^{m+1}$,
there exist $t_1, \dots ,t_m$ in $R$ such that
$(r_1+t_1r_{m+1}, \dots ,r_m+t_mr_{m+1}) \in R^m$ is unimodular. We say
that $R$ has {\it stable rank} $m$, we denote it with ${\rm sr}(R)=m$,
if $m$ is
the least number such that $({\rm S}_m)$ holds. If such a number does not
exist we say that ${\rm sr}(R)=\infty$.

An $n\times k$-matrix $B$ with $n<k$ is called unimodular if $B$ has a right
inverse. If $B$ is an $n\times k$-matrix and $C \in \GL(k, R)$, then
$B$ is unimodular if and only if  $CB$ is unimodular.
A matrix of the form
$\left(\begin{array}{cc}
1 & 0      \\
u &  B
\end{array} \right)$,
where $u$ is a column vector with coordinates in $R$, is unimodular
if and only if the matrix $B$ is unimodular.

We say that the ring $R$ satisfies the
{\it stable range condition} $({\rm S}_n^k)$
if for every $n\times (n+k)$-matrix $B$, there exists a vector
$r=(r_1, \dots, r_{n+k-1})$ such that
$B\left(\begin{array}{cc}
1 &    r     \\
0  &  I_{n+k-1}
\end{array} \right)
=\left(\begin{array}{cc}
u & B'
\end{array} \right)$,
where the $n \times (n+k-1)$-matrix $B'$ is unimodular and $u$ is the
first column of the matrix $B$.

\begin{thm}[Vaserstein]\label{vas3}
For every $k \ge 1$ and $n \ge 1$, a ring $R$ satisfies $({\rm S}_k)$
if and only if it satisfies $({\rm S}_n^k)$.
\end{thm}
\begin{proof}
The definition of $({\rm S}_n^k)$ and this theorem are just the transposed
version of theorem \cite[Thm.~3$'$]{vas4} of Vaserstein.
\end{proof}

Let $\u (R^n)$ denote the
subposet of $\o(R^n)$ consisting of unimodular sequences. Recall that
a sequence of vectors $v_1, \dots, v_k$ in $R^n$ is
called unimodular when $v_1, \dots, v_k$ is basis of a free direct
summand of $R^n$. Note that if $(v_1, \dots, v_k) \in \o(R^n)$ and if
$n \le m$, it is the
same to say that $(v_1, \dots, v_k)$ is unimodular as a sequence of vectors
in $R^n$ or as a sequence of vectors in $R^m$. We call an element
$(v_1, \dots, v_k)$ of $\u(R^n)$ a $k$-{\it frame}.

\begin{thm}[Van der Kallen]\label{kal5}
Let $R$ be a ring with ${\rm sr}(R) < \infty$ and $n \le m$. Then
\par {\rm (i)} $\o(R^n) \cap \u(R^m)$ is $(n-{\rm sr}(R)-1)$-connected.
\par {\rm (ii)} $\o(R^n) \cap \u(R^m)_v$ is
$(n-{\rm sr}(R)-|v|-1)$-connected for all $v \in \u(R^m)$.
\end{thm}
\begin{proof}
See \cite[Thm. 2.6]{kall}.
\end{proof}

{}From now on let $R$ be a commutative ring.
Let $e_{i, j}(r)$ be the $2n \times 2n$-matrix with $r \in R$ in
the $(i, j)$ place and zero elsewhere.
Consider $Q= \s_{i=1}^n(e_{2i-1, 2i}(1)-e_{2i, 2i-1}(1)) \in \GL(2n, R)$.
By definition the symplectic group is the group
\begin{gather*}
\Spp(2n, R):=\{ A \in \GL(2n, R) : {}^t\!AQA=Q\}.
\end{gather*}
Let $\si$ be the permutation of the set of natural numbers given by
$\si(2i)=2i-1$
and  $\si(2i-1)=2i$. For $1 \le i, j \le 2n$, $i \neq j$, and every
$r \in R$ define
\begin{gather*}
E_{i, j}(r)=\begin{cases}I_{2n}+ e_{i, j}(r) & \text{if $i=\si(j)$}\\
I_{2n}+ e_{i, j}(r)-(-1)^{i+j}e_{\si(i), \si(j)}(r)& \text{if $i \neq\si(j)$
and $i<j$} \end{cases}
\end{gather*}
where $I_{2n}$ is the identity element of $\GL(2n, R)$.
It is easy to see that $E_{i, j}(r) \in \Spp(2n, R)$, for every $r \in R$.
Let  $\ESp(2n, R)$ be the group generated by $E_{i, j}(r)$, $r \in R$.
We call it  {\it elementary
symplectic group}. Define the
bilinear map $h: R^{2n} \times R^{2n} \arr R$ by
$h(x, y)={}^t\!xQy= \s_{i=1}^n(x_{2i-1}y_{2i}-y_{2i-1}x_{2i})$ where
$x=(x_1, \dots, x_{2n})$ and $y=(y_1, \dots, y_{2n})$. Clearly $h(x, x)=0$
for every $x \in R^{2n}$. We say that a subset $S$ of $R^{2n}$ is isotropic
if for every $x, y \in S$, $h(x, y)=0$. So every
element of $\rn$ is isotropic. 
If $h(x, y)=0$, then we say that $x$ is perpendicular to $y$.
We denote by $\lan S \ran$ the submodule of $\rn$ generated by $S$,
and by $\lan S \ran^\perp$ the submodule consisting of all the
elements of $\rn$ which are perpendicular to all the elements of $S$.

Let $\iur$ be the set of sequences $(x_1, \dots, x_k)$, $x_i \in R^{2n}$,
such that $x_1, \dots, x_k$ form a basis for an isotropic direct summand
of $R^{2n}$. Let
$\hur$ be the set of sequences $((x_1, y_1), \dots, (x_k, y_k))$
such that $(x_1, \dots, x_k)$, $(y_1, \dots, y_k) \in \iur$,
$h(x_i, y_j)=\delta_{i,j}$, where $\delta_{i,j}$ is the Kronecker delta.
We call $\iur$ and $\hur$ the
poset of isotropic unimodular sequences and the poset of
hyperbolic unimodular sequences, respectively.
For $1 \le k \le n$, let
$\iu(\rn, k)$ and $\hu(\rn, k)$ be the set of all elements of length $k$
of $\iur$ and $\hur$ respectively.
We call the elements of $\iu(\rn, k)$ and $\hu(\rn, k)$ the
isotropic $k$-frames and the hyperbolic $k$-frames, respectively.
Define the poset $\mur$ as the set of
$((x_1, y_1), \dots, (x_k, y_k)) \in \o(\rn \times \rn)$ such that,
{\rm (i)} $(x_1, \dots, x_k) \in \iur$,
{\rm (ii)} for each $i$, either $y_i=0$ or $(x_j, y_i)=\delta_{ji}$,
{\rm (iii)} $\lan y_1, \dots, y_k \ran$ is isotropic.
We identify
$\iur$ with $\mur \cap \o(\rn \times \{0\})$ and
$\hur$ with $\mur \cap \o(\rn \times (\rn \bs \{0\}))$.

\begin{lem}\label{vas0}
Let $R$ be a ring with $\sr < \infty$. If $n \ge \sr + k$ then
$\ESp(2n, R)$ acts transitively on $\iu(\rn, k)$ and $\hu(\rn, k)$.
\end{lem}
\begin{proof}
The proof is by induction on $k$. If $k=1$, the fact that $\ESp(2n, R)$ acts
transitively on  $\iu(\rn, 1)$, is due to Vaserstein \cite[Thm. 2.]{vas3}.
The transitivity of action  on $\hu(\rn, 1)$ is also
well known and easy. In fact if $(x, y)$ is a hyperbolic pair then there is
an $E \in \ESp(2n, R)$ such that $Ex=e_1$. Thus
$Ey=(r_1, 1, r_3, \dots , r_{2n})$. Now if
$E'=E_{1, 2}(-r_1) \prod_{i=3}^{2n}E_{i, 2}(-r_i)$, then $E' E x=e_1$ and
$E' E y=e_2$. This shows that $\ESp(2n, R)$ acts transitively on  $\hu(\rn,
1)$.
The rest is an easy induction and the fact that for every isotropic $k$-frame
$(x_1, \dots, x_k)$ there is an  isotropic $k$-frame
$(y_1, \dots, y_k)$ such that $((x_1, y_1), \dots, (x_k, y_k))$ is a
hyperbolic $k$-frame.
\end{proof}

\begin{lem}\label{charn}
Let $R$ be a ring with $\sr < \infty$, and let $n \ge \sr + k$.
Let $((x_1, y_1), \dots, (x_k, y_k))$ $\in \hur$, $(x_1, \dots, x_k)$ $ \in
\iur$
and $V=\lan x_1, \dots, x_k\ran$. Then
\par {\rm(i)} $\iur_{(x_1, \dots, x_k)} \simeq \iu(R^{2(n-k)})\lan V \ran$,
\par {\rm(ii)} $\hur \cap \mur_{((x_1, 0), \dots, (x_k, 0))}
\simeq \hur_{((x_1, y_1), \dots, (x_k, y_k))}\lan V \times V \ran$,
\par {\rm(iii)} $\hur_{((x_1, y_1), \dots, (x_k, y_k))} \simeq
\hu(R^{2(n-k)})$.
\end{lem}
\begin{proof}
See \cite{ch}, the proof of lemma 3.4 and the proof of Thm. 3.2.
\end{proof}

\begin{lem}\label{b-w0}
Let $R$ be a ring with $\sr < \infty$. Assume $n \ge \sr + k$ and
$(v_1,\dots, v_k) \in \u(\rn)$. Then there is a hyperbolic basis
$\{x_1, y_1, \dots x_n, y_n\}$ of $\rn$ such that
$v_1, \dots, v_k \in \lan x_1, y_1, \dots x_{k-1}, y_{k-1}, x_k\ran$.
\end{lem}
\begin{proof}
By lemma \ref{vas0}, we can assume that $v_1=e_1$. 
The proof is by induction on $k$.
If $k=1$, everything is trivial.
Consider the $k \times 2n$-matrix
$A$, whose $i$-th row is the vector $v_i$. Let $B$ be the
$(k-1) \times (2n-1)$-matrix obtained from $A$
by deleting the first column and the
first row. Since $A$ is unimodular, $B$ is
unimodular too. By theorem \ref{vas3}, there exist a vector
$s=(s_3, \dots, s_{2n})$ such that
$B\left(\begin{array}{cc}
1 &    s     \\
0  &  I_{2n-2}
\end{array} \right)
=\left(\begin{array}{cc}
u & C
\end{array} \right)$,
where the $(k-1)\times (2n-2)$-matrix $C$ is unimodular and $u$ is the
first column of the matrix $B$. Now let
$E=\prod_{i=3}^{2n}E_{2, i}(s_i) \in \ESp(2n, R)$. An easy computation
shows that
\begin{gather*}
AE=\left(\begin{array}{cccc}
1 & 0 & \dots &  0   \\
\ast & * &   &  \\
\vdots &\vdots &  & C \\
\ast &  * &   &
\end{array} \right).
\end{gather*}
Denote the rows of $C$ by $w_1, \dots, w_{k-1} \in R^{2n-2}$. Since
$n-1 \ge \sr + (k-1)$, by induction hypothesis there exist a hyperbolic basis
$\{x_2, y_2, \dots, x_n, y_n\}$ of $R^{2n-2}$ such that
$w_i \in \langle x_2, y_2, \dots x_{k-1}, y_{k-1}, x_{k}\rangle$.
If we consider $R^{2n-2}$ as a subspace of $\rn$ generated by
$e_3, \dots, e_{2n}$ and if $x_1=e_1$ and
$y_1=e_2$, then $v_i \in \langle x_1, y_1, \dots x_{k-1}, y_{k-1},
x_{k}\rangle$. This completes the proof.
\end{proof}

\begin{lem}\label{u-i}
Let $R$ be a ring with $\sr < \infty$.
Let  $(v_1, \dots, v_k) \in \u(\rn)$.
If $n \ge \sr + k$ then
\par {\rm (i)} $\u(\rn)_{(v_1, \dots, v_k)} \cap
\o(\langle v_1, \dots, v_k\rangle^\perp)$ is $(2n-\sr-2k-1)$-connected.
\par {\rm (ii)} ${\u(\rn)}_{(v_1, \dots, v_k, w_1, \dots, w_r)} \cap
\o(\langle v_1, \dots, v_k \rangle^\perp)$ is $(2n-\sr-2k-r-1)$-connected,
for every $(w_1, \dots, w_r) \in {\u(\rn)}_{
(v_1,\dots, v_k)}$.
\end{lem}
\begin{proof}
It is well known and easier than lemma \ref{vas0} that $\GL(2n,R)$
acts transitively on $k$-frames. Therefore let us choose a basis
$b_1,\ldots,b_{2n}$ of $\rn$ so that $(b_{2n-k+1},\ldots, b_{2n})=
(v_1, \dots, v_k)$. Let $c_1,\ldots , c_{2n}$ be the dual basis:
$h(b_i, c_j)=\delta_{i,j}$, where $\delta_{i,j}$ is the Kronecker delta.
Then $\langle v_1, \dots, v_k\rangle^\perp=c_1R+\dots +c_{2n-k}R$,
so we may identify $\o(\langle v_1, \dots, v_k\rangle^\perp)$
with $\o(R^{2n-k})$ after a change of basis.  Now by
theorem \ref{kal5}, the poset $\u(\rn)_{(v_1, \dots, v_k)} \cap
\o(R^{2n-k})$
is $(2n-\sr-2k-1)$-connected.
The proof of (ii) is similar to the proof of (i).
\end{proof}

For a real number $l$, by $\floor [l]$ we mean the largest integer $n$ with
$n\leq l$.

\begin{thm}\label{b-w1}
The poset $\iur$ is
$\floor [\frac{n-\sr-2}{2}]$-connected and $\iur_x$ is
$\floor [\frac{n-\sr-|x|-2}{2}]$-connected for every
$x \in \iur$.
\end{thm}
\begin{proof}
If $n\leq\sr$, the result is clear,
so let $n>\sr$. Let
$X_v=\iur \cap \u(\rn)_v \cap \o(\langle v \rangle^\perp)$,
for every $v \in \u(\rn)$, and put $X:= \bcu_{v \in F}X_v$ where $F=\u(\rn)$.
It follows from lemma \ref{vas0}
that $\iur_{\le n-\sr} \se  X$.
So to treat $\iur$, it is enough to prove that
$X$ is $\floor [\frac{n-\sr-2}{2}]$-connected.
First we prove that  $X_v$ is
$\floor [\frac{n-\sr-|v|-2}{2}]$-connected for every $v \in F$.
The proof is by descending
induction on $|v|$.
If $|v| > n-\sr$, then $\floor [\frac{n-\sr-|v|-2}{2}]<-1$.
In this case there
is nothing to prove. If $n-\sr-1\le|v| \le n-\sr$, then
$\floor [\frac{n-\sr-|v|-2}{2}]= -1$,
so we must prove that $X_v$ is nonempty. This follows from lemma \ref{b-w0}.
Now assume $|v|\leq n-\sr-2$ and assume
by induction that $X_w$ is $\floor [\frac{n-\sr-|w|-2}{2}]$-connected
for  every $w$, with $|w| >|v|$.
Let $l=\floor [\frac{n-\sr-|v|-2}{2}]$, and observe that $ n-|v|-\sr\geq l+2$.
Put $T_w=\iur \cap \u(\rn)_{wv} \cap \o(\lan wv \ran^\perp)$
where $w \in G_v=\u(\rn)_v\cap \o(\langle v \rangle^\perp)$
and put $T:=\bcu_{w \in G_v}T_w$ . It follows by lemma \ref{b-w0}
that $(X_v)_{\le n-|v|-\sr} \se T$.
So it is enough to prove that $T$ is $l$-connected.
The poset $G_v$ is $l$-connected by lemma \ref{u-i}.
By induction,
$T_w$ is $\floor [\frac{n-\sr-|v|-|w|-2}{2}]$-connected.
But ${\rm min}\{l-1, l-|w|+1\} \le \floor [\frac{n-\sr-|v|-|w|-2}{2}]$,
so $T_w$ is
${\rm min}\{l-1, l-|w|+1\}$-connected. For every $ y \in T$,
$\a_y= \{w \in G_v: y \in T_w\}$ is isomorphic to
$\u(\rn)_{vy} \cap \o(\langle
vy \rangle^\perp)$ so by lemma \ref{u-i},
it is $(l-|y|+1)$-connected.
Let $w\in G_v$ with $|w|=1$.
For every $z\in T_w$ we have $wz\in X_v$,
so $T_w$ is contained in a cone,  call it
$C_w$, inside $X_v$. Put $C(T_w)=T_w \cup (C_w)_{\le n-|v|-\sr}$.
Thus $C(T_w) \se T$. The poset $C(T_w)$ is $l$-connected because
$C(T_w)_{\le n-|v|-\sr} = (C_w)_{\le n-|v|-\sr}$.
Now by theorems \ref{kal5} and \ref{surj}, $T$ is $l$-connected.
In other words, we have now shown that
$X_v$ is $\floor [\frac{n-\sr-|v|-2}{2}]$-connected.
By knowing this one can prove, in a similar way, that $X$ is
$\floor [\frac{n-\sr-2}{2}]$-connected. (Just pretend that $|v|=0$.)

Now consider the poset $\iur_x$ for an $x=(x_1, \dots, x_k) \in \iur$.
The proof is by induction on $n$. If $n=1$, everything is easy.
Similarly, we may assume $n-\sr-|x|\geq 0$.
Let $l=\floor [\frac{n-\sr-|x|-2}{2}]$. By lemma \ref{charn},
$\iur_x \simeq \iu(R^{2(n-|x|)})\lan V \ran$, where
$V=\lan x_1, \dots, x_k \ran$. In the above we proved that $\iu(R^{2(n-|x|)})$
is $l$-connected and by induction, the poset
$\iu(R^{2(n-|x|)})_y$ is $\floor [\frac{n-|x|-\sr-|y|-2}{2}]$-connected
for every $y \in \iu(R^{2(n-|x|)})$.
But $l-|y| \le \floor [\frac{n-|x|-\sr-|y|-2}{2}]$. So
$\iu(R^{2(n-|x|)})\lan V \ran$ is $l$-connected  by
lemma \ref{maazen5}. Therefore $\iur_x$
is $l$-connected.
\end{proof}

\begin{thm}\label{b-w2}
The poset $\hur$ is $\floor [\frac{n-\sr-3}{2}]$-connected and
$\hur_x$ is $\floor [\frac{n-\sr-|x|-3}{2}]$-connected for every
$x \in \hur$.
\end{thm}
\begin{proof}
The proof is by induction on $n$. If $n=1$, then everything is trivial.
Let $F=\iur$ and $X_v=\hur \cap \mur_v$,
for every $v \in F$. Put  $X:=\bcu_{v \in F} X_v$.
It follows from lemma \ref{vas0} that $\hur_{\le n-\sr} \se X$.
Thus to treat $\hur$, it is enough to prove that $X$ is
$\floor [\frac{n-\sr-3}{2}]$-connected, and we may assume $n\geq\sr+1$.
Take
$l=\floor [\frac{n-\sr-3}{2}]$ and $V=\lan v_1, \dots, v_k \ran$, where
$v=(v_1, \dots, v_k)$. By lemma \ref{charn}, there is an isomorphism
$X_v \simeq \hu(R^{2(n-|v|)}) \lan V \times V \ran$, if $n\geq\sr+|v|$.
By induction
$\hu(R^{2(n-|v|)})$ is $\floor [\frac{n-|v|-\sr-3}{2}]$-connected and
again by induction
$\hu(R^{2(n-|v|)})_y$ is
$\floor [\frac{n-|v|-\sr-|y|-3}{2}]$-connected for every
$y \in \hu(R^{2(n-|v|)})$. So by  lemma \ref{maazen5}, $X_v$ is
$\floor [\frac{n-|v|-\sr-3}{2}]$-connected.
Thus the poset $X_v$ is
${\rm min}\{l-1, l-|v|+1\}$-connected.
Let
$x=((x_1, y_1), \dots, (x_k, y_k))$. It is easy to see that
$\a_x=\{ v \in F: x \in X_v\} \simeq \iur_{(x_1, \dots, x_k)}$.
By the above theorem \ref{b-w1}, $\a_x$ is
$\floor [\frac{n-\sr-k-2}{2}]$-connected. But
$l-|x|+1 \le \floor [\frac{n-\sr-k-2}{2}]$, so $\a_x$ is $(l-|x|+1)$-connected.
Let $v=(v_1) \in F$, $|v|=1$, and let
$D_v:= \hur_{(v_1, w_1)}\simeq \hu(R^{2(n-1)})$ where $w_1 \in \rn$ is a
hyperbolic
dual of $v_1 \in \rn$.
Then $D_v \se X_v$ and $D_v$ is contained in a
cone, call it $C_v$, inside $\hur$. Take
$C(D_v):=D_v \cup (C_v)_{\le n-\sr}$.
By induction $D_v$ is $\floor [\frac{n-1-\sr-3}{2}]$-connected and so
$(l-1)$-connected.
Let $Y_v= X_v \cup C(D_v)$. By the Mayer-Vietoris
theorem and the fact that $C(D_v)$ is $l$-connected,
we get the exact sequence
\[
\H_l(D_v, \z) \overset{(i_v)_\ast}{\arr} \H_l(X_v, \z)
\arr \H_l(Y_v, \z) \arr 0.
\]
where $i_v: D_v \arr X_v$ is the inclusion.
By induction $(D_v)_w$ is $\floor [\frac{n-1-\sr-|w|-3}{2}]$-connected
and so
$(l-|w|)$-connected, for $w\in D_v$.
By lemma \ref{maazen5}(i)
and lemma \ref{charn},
$(i_v)_\ast$ is an isomorphism,
and by exactness of the above sequence we get
$\H_l(Y_v, \z)=0$. If $l \ge 1$ by the Van Kampen theorem
$\pi_1(Y_v, x) \simeq \pi_1(X_v, x)/N$ where $x \in D_v$ and $N$ is the
normal subgroup generated by the image of the map
$(i_v)_\ast: \pi_1(D_v, x) \arr \pi_1(X_v, x)$. Now by lemma
\ref{maazen5}(ii), $\pi_1(Y_v, x)$ is trivial.
Thus by the Hurewicz theorem \ref{hur}, $Y_v$ is $l$-connected.
By having all this we can apply  theorem
\ref{surj} and so $X$ is $l$-connected.
The fact that $\hur_x$ is
$\floor [\frac{n-\sr-|x|-3}{2}]$-connected follows from the above
and lemma \ref{charn}.
\end{proof}

\begin{rem}
Charney in \cite[2.10]{ch}, conjectured that $\iur$, with respect
to certain bilinear forms, is highly connected. Our theorem \ref{b-w1} proves
this conjecture for the above bilinear form.
Also, by assuming the high connectivity of the $\iur$, she
proved that $\hur$ is highly connected. Our proof is different 
and relies on our theory, but we use ideas from her paper,
such as the lemma \ref{charn} and her  lemma \ref{maazen5},
which is a modified version of work of
Maazen.
\end{rem}

\begin{rem}
One should notice that in the above, we used heavily the fact that
every unimodular vector (in fact any vector) in $\rn$ is isotropic.
So our arguments don't work for orthogonal groups. One needs to
do more. Consider the bilinear form $h': \rn \times \rn \arr R$
defined by $h'(x, y)= \s_{i=1}^n(x_{2i-1}y_{2i}+y_{2i-1}x_{2i})$. Let
$T=\{ x \in \rn: h'(x, x)=0\}$ and $\u'(\rn)=\u(\rn) \cap \o(T)$.
If $\u'(\rn)$ is highly connected then we have similar results
for the poset $\iur$
arising from the form $h'$. So we formulate
the following question:\\
Question. Is $\u'(\rn)$ is highly connected for example $(n-\sr-1)$-connected?
\end{rem}

\section{ Homology stability}

{}From theorem \ref{b-w2} one can get the homology stability of
symplectic groups as Charney proved in \cite[Sec. 4]{ch}. Here we only
formulate the theorem and for the proof we refer to Charney's paper.

\begin{thm}
Let $R$ be a commutative ring with $ \sr < \infty$. Then for every abelian
group $L$ the homomorphism
$ {\psi_n}_\ast : H_i(\Spp(2n, R), L) \arr H_i(\Spp(2n+2, R), L)$
is surjective for $n \ge 2i+\sr+2$ and bijective for $n \ge 2i+\sr+3$, where
$ {\psi_n} : \Spp(2n, R) \arr \Spp(2n+2, R)$, $A \mt
\left(\begin{array}{cc}
A & 0      \\
0 &  I_2
\end{array} \right)$.
\end{thm}
\begin{proof}
See \cite[Sec. 4]{ch}.
\end{proof}

\begin{rem}
To prove homology stability of this type one only needs high
acyclicity of the corresponding poset. But usually this type of posets are
also highly connected and it looks that it is a tradition to give a
proof of high connectivity of the posets. In this paper
we follow the tradition. In particular we
wished to confirm  the conjecture of Charney \cite[2.10]{ch}
for the particular bilinear form that we considered.
\end{rem}

\

 e-mail:\quad \texttt{mirzaii@math.uu.nl\qquad vdkallen@math.uu.nl}
\end{document}